\documentclass[11pt,a4paper,equation]{article}
\usepackage[greek,american]{babel}
\usepackage[iso-8859-7]{inputenc}
\usepackage{amssymb,amsfonts}
\usepackage[dvips]{graphicx}
\usepackage{amsmath}
\usepackage{epsfig}
\usepackage{latexsym}
\usepackage{makeidx}
\usepackage{pst-math,pst-xkey}

\numberwithin{equation}{section}
\makeindex
\begin{document}
\date{}
\author{\textbf{Vassilis G. Papanicolaou$^1$ and Aristides V. Doumas$^2$}
\\\\
Department of Mathematics
\\
National Technical University of Athens,
\\
Zografou Campus, 157 80, Athens, GREECE
\\\\
$^1${\tt papanico@math.ntua.gr} \qquad  $^2${\tt adou@math.ntua.gr}}
\title{On an Old Question of Erd\H{o}s and R\'{e}nyi Arising in
\\
the Delay Analysis of Broadcast Channels}
\maketitle
\begin{abstract}
Consider a broadcast channel with $n$ users, where different users receive different messages, and suppose that each user has to receive $m$
packets. A quantity of interest here, introduced by Sharif and Hassibi (2006-7) \cite{Sh}, \cite{S-H}, is the (\emph{packet}) \emph{delay}
$D_{m,n}$, namely the number of channel uses required to guarantee that all users will receive $m$ packets. For the case of a \emph{homogeneous}
network, where in each channel use the transmitter chooses a user at random, i.e. with probability $1/n$, and sends him$/$her a packet,
the same quantity $D_{m,n}$ had already appeared in the \emph{coupon collector} context, in the works of Newman and Shepp
(1960) \cite{N-S} and of Erd\H{o}s and R\'{e}nyi (1961) \cite{E-R}.

A problem of particular interest in wireless communications, related to the delay $D_{m,n}$, is to determine its behavior as $n$ and $m$
grow large. Regarding this problem, Sharif and Hassibi \cite{Sh}, \cite{S-H} managed to calculated the asymptotics of the mean value
$\mathbb{E}[D_{m,n}]$, as $n \to \infty$, for the cases (a) $m = \ln n$ and (b) $m = (\ln n)^{\rho}$, $\rho > 1$.
It is remarkable that Erd\H{o}s and R\'{e}nyi \cite{E-R} had, also,
raised the question of the determination of the asymptotic profile of $D_{m,n}$ for large $m$ and $n$ (in 1961). And in the 1970's the limiting distribution of $D_{m,n}$ for large $m$ and $n$ was determined (in great generality) by Ivchenko \cite{I1} and \cite{I2}.

In this article we determine the asymptotics of the moments of $D_{m,n}$ for large $m$ and $n$. We also derive its limiting distribution in the
``supercritical case" where $m$ grows faster than $\ln n$ and in the ``critical case" $m \sim \beta \ln n$, by an approach which is different from the one used by Ivchenko.
\end{abstract}

\textbf{2010 AMS Mathematics Classification.} 60F05; 60F99.

\section{Introduction}
In their works \cite{Sh}, \cite{S-H} M. Sharif and B. Hassibi consider a single-antenna broadcast fading channel with $n$ users, where the
transmission is packet-based. They define the (\emph{packet}) \emph{delay} $D_{m,n}$ as the minimum number of channel uses that guarantees all $n$
users to successfully receive $m$ packets (clearly, $D_{m,n} \geq mn$). Sharif and Hassibi consider an ``opportunistic" scheduling, where in each
channel use the transmitter sends the packet to the user with the best channel conditions, i.e. the highest signal-to-noise ratio (SNR). As it turns
out \cite{S-H}, the opportunistic scheduling maximizes the sum rate (or \emph{throughput}), namely the rate of successful message delivery of the
broadcast channel.

The case where all users have the same SNR is referred as the \emph{homogeneous} network case. Here, in each channel use the transmitter chooses the
$j$-th user, $j = 1, 2, \ldots, n$, with probability $1/n$. It follows that $D_{m,n}$ becomes a classical quantity related to \emph{urn problems} and, in particular, to the \emph{coupon collector's problem} (\emph{CCP}): Suppose $n$ equally likely coupons are sampled independently with replacement.
Then, $D_{m,n}$ is the number of trials needed until all $n$ coupons are detected at least $m$ times.
If the coupon probabilities are not equal, then the CCP-quantity $D_{m,n}$ corresponds to the delay of a \emph{heterogeneous} network.
In the present paper, however, we will only consider the homogeneous network case.

In the CCP context, the random variable $D_{m,n}$ first appeared in the works \cite{N-S} of D.J. Newman and L. Shepp ($1960$), and
\cite{E-R} of P. Erd\H{o}s and A. R\'{e}nyi ($1961$). Both works focused on the asymptotic behavior of $D_{m,n}$ as $n \to \infty$, while $m$ stays
fixed. Newman and Shepp \cite{N-S} obtained that
\begin{equation}
\mathbb{E}\left[D_{m,n}\right] = n \ln n + \left(m-1\right) n \ln \ln n + n C_m + o(n)
\label{1}
\end{equation}
as $n\rightarrow \infty$, where $C_m$ is a constant depending on $m$. Roughly speaking, formula (\ref{1}) tells us that, on the average, the detection of all $n$ coupons at least once ``costs" $n \ln n + O(n)$, while each additional detection (of all coupons) raises the cost by
$n \ln \ln n  + O(n)$.

Soonafter, Erd\H{o}s and R\'{e}nyi \cite{E-R} went a step further and determined the limit distribution of $D_{m,n}$, as well as the exact value of
the constant $C_{m}$. They proved that
 for every real $y$ one has
\begin{equation}
\lim_{n \rightarrow \infty} \mathbb{P}\left\{\frac{D_{m,n} - n \ln n - (m - 1) n \ln \ln n}{n} \leq y\right\}
= \exp \left(-\frac{e^{-y}}{(m-1)!}\right)
\label{3}
\end{equation}
or, equivalently,
\begin{equation}
\lim_{n \rightarrow \infty} \mathbb{P}\left\{\frac{D_{m,n}}{n} - \ln n - (m - 1) \ln \ln n + \ln (m-1)!\leq y\right\}
= e^{-e^{-y}}.
\label{333}
\end{equation}
Noticing that in the right-hand side of \eqref{333} we have the standard Gumbel distribution function, whose expectation is
$\gamma = 0.5772\cdots$ (the Euler-Mascheroni constant), it is not hard to justify that the constant in \eqref{1} must be
\begin{equation}
C_m = \gamma - \ln \left(m-1\right)!.
\label{2}
\end{equation}


At the end of their paper, Erd\H{o}s and R\'{e}nyi \cite{E-R} have included the following comment:

\smallskip

\textit{``It is an interesting problem to investigate the limiting distribution of
$\nu_m(n)$ when $m$ increases with $n$, but we can not go into this question here."}

\smallskip

The quantity denoted by $\nu_m(n)$ in \cite{E-R} is nothing but $D_{m,n}$.

Formulas \eqref{1}, \eqref{2}, and \eqref{333} hint that the case where $m$ grows like $\ln n$ is expected to be ``critical"
(as opposed to the ``subcritical" case where $m$ grows slower than $\ln n$) in the sense that it is the smallest growth of $m$ (with respect to $n$)
which seems to affect the leading asymptotic behaviors of $\mathbb{E}\left[D_{m,n}\right]$ and $D_{m,n}$. The above guess was confirmed by
G.I. Ivchenko, who, actually, managed to give a complete answer to the question of Erd\H{o}s and R\'{e}nyi in his works \cite{I1} and \cite{I2}
(in the 1970's).

Sharif and Hassibi \cite{Sh}, \cite{S-H} were also interested in the behavior of $\mathbb{E}\left[D_{m,n}\right]$ as both $m$ and $n$ grow to
infinity (although they did not seem to be aware neither of the paper of Erd\H{o}s and R\'{e}nyi), nor of the works of Ivchenko.
Sharif and Hassibi managed to show \cite{Sh}, \cite{S-H} that:

(i) If $m = \ln n$ and $n \to \infty$, then
\begin{equation}
\mathbb{E}\left[D_{m,n}\right] = \alpha \, n \ln n + O(n \ln\ln n),
\label{SH1}
\end{equation}
where $\alpha$ is the (unique) solution of the equation $\alpha - \ln \alpha = 2$ in the interval $(1, \infty)$ ($\alpha \approx 3.146$).

(ii) If $m = (\ln n)^{\rho}$, where $\rho > 1$ is fixed and $n \to \infty$, then
\begin{equation}
\mathbb{E}\left[D_{m,n}\right] = nm + o(mn).
\label{SH2}
\end{equation}

They also showed that if $m \to \infty$, while $n$ stays fixed, then $\mathbb{E}[D_{m,n}] = nm + o(m)$, but this fact had beed
already known to Newman and Shepp \cite{N-S} (Sharif and Hassibi, though, seem to be aware of the fact that \eqref{SH2}
remains valid whenever $m$ grows faster than $\ln n)$.

The reader may have noticed that, since $m$ and $n$ are integers, the above mentioned equalities $m = \ln n$ and $m = (\ln n)^{\rho}$, strictly
speaking, cannot be satisfied; they only make sense asymptotically (e.g., $m = \ln n$ can be interpreted as $m = \ln n + O(1)$).

In Section 2 of the present paper we calculate the (leading) asymptotics of the moments of $D_{m,n}$ in the cases (i) $m \gg \ln n$ and
(ii) $m \sim \beta\ln n$, where $\beta > 0$ (as usual, the notation $A(n) \gg B(n)$ means that $B(n) / A(n) \to 0$ as $n \to \infty$).
Then, in Section 3, under some mild restrictions on the growth of $m$ in the aforementioned cases (i) and (ii), namely (i$\,\acute{}\,$) $m \gg \ln^3 n$ and
(ii$\,\acute{}\,$) $m = \beta\ln n + o\big(\sqrt{\ln n}\,\big)$, we determine the limiting distribution of $D_{m,n}$ by a method which is
different from that of Ivchenko. As we will see, the quantity $D_{m,n}$, appropriately normalized, converges in distribution to a
Gumbel random variable. This may not sound surprising, but the challenging part is to obtain the right normalization of $D_{m,n}$.

\subsection{Analyzing $D_{m,n}$ via Poissonization}

Suppose that, for $j = 1, \dots, n$, we denote by $X_j$ the number of trials needed in order to detect the $j$-th coupon $m$ times.
Then, it is clear that $X_j$ is a negative binomial random variable, with parameters $m$ and $1/n$, and
\begin{equation*}
D_{m,n} = \max_{1 \leq j \leq n} X_j.
\end{equation*}
However, the above formula for $D_{m,n}$ is not very convenient, since the $X_j$'s are not independent. Fortunately, there is a clever
``Poissonization technique" (probably due to L. Holst \cite{Ho} --- see also \cite{Ro2}) from which we can get more insight about $D_{m,n}$.

Let $N(t)$, $t \geq 0$, be a Poisson process with rate $1$. We imagine that each Poisson event associated to this process is a sampled
coupon, so that $N(t)$ is the number of sampled coupons at time $t$.
Next, for $j = 1, \dots, n$, let $N_j(t)$ be the number of detections of the $j$-th coupon at time $t$.
Then, the processes $N_j(t)$, $j = 1, \dots, n$, are independent Poisson processes with rates $1/n$
\cite{Ro2} and it is clear that $N(t) = N_1(t) + \cdots + N_n(t)$. If $T_j$, $j = 1, \dots, n$, denotes the time of the $m$-th
event of the process $N_j$, then $T_1, \dots, T_n$ are independent (being associated to independent processes) and
\begin{equation}
\Delta_{m,n} := \max_{1 \leq j \leq n} T_j
\label{P1}
\end{equation}
is the time when all different coupons have been detected at least $m$ times.

Now, for each $j = 1, \dots, n$, the random variable $T_j$ (being a sum of $m$ independent exponential variables with parameter $1/n$) is Erlang with parameters $m$ and $1/n$, hence
\begin{equation}
\mathbb{P}\{T_j > t\} = S_m(t/n) e^{-t/n},
\qquad
t \geq 0,
\label{P2}
\end{equation}
where
$S_{m}(y)$ denotes the $m$-th partial sum of the Taylor-Maclaurin series of $e^{y}$, namely
\begin{equation}
S_{m}(y) := 1+y+\frac{y^{2}}{2!}+\cdots+\frac{y^{m-1}}{\left(m-1\right)!}=\sum_{l=0}^{m-1}\frac{y^l}{l!}
\label{7}
\end{equation}

It follows from \eqref{P1}, \eqref{P2}, and the independence of the $T_j$'s that the distribution function of $\Delta_{m,n}$ is
\begin{equation}
F_{\Delta}(t) = \mathbb{P}\{\Delta_{m,n} \leq t\} = \mathbb{P}\{T_1 \leq t\}^n = \left[ 1 - S_m(t/n) e^{-t/n} \right]^n,
\qquad
t \geq 0.
\label{P3}
\end{equation}
It remains to relate $\Delta_{m,n}$ to $D_{m,n}$. Clearly,
\begin{equation}
\Delta_{m,n} = \sum_{k=1}^{D_{m,n}} U_k,
\label{P4}
\end{equation}
where $U_1, U_2, \dots$ are the interarrival times of $N(t)$. Since $N(t)$ is a Poisson process of rate $1$, the $U_j$'s are independent exponential
random variables with parameter $1$. Furthermore, it is clear that $D_{m,n}$ is independent of the $U_j$'s \cite{Ro2}.

One consequence of formula \eqref{P4} is that, given $D_{m,n}$, the variable $\Delta_{m,n}$ is Erlang with parameters $D_{m,n}$ and $1$, i.e.
\begin{equation}
\mathbb{P}\left\{\Delta_{m,n} > t \, | \, D_{m,n}\right\} = S_{D_{m,n}}(t) e^{-t},
\qquad
t \geq 0.
\label{A1}
\end{equation}
In other words, the conditional probability density of $\Delta_{m,n}$ given $D_{m,n}$ is
\begin{equation}
f_{\Delta | D}(t) = \mathbb{P}\left\{\Delta_{m,n} \in dt \, | \, D_{m,n}\right\} = \frac{t^{D_{m,n}-1}}{\left(D_{m,n}-1\right)!} \, e^{-t},
\qquad
t \geq 0.
\label{A2}
\end{equation}
We can take expectations in \eqref{A1} and obtain
\begin{equation}
 \mathbb{E}\left[S_{D_{m,n}}(t)\right] e^{-t} = \mathbb{P}\left\{\Delta_{m,n} > t\right\} = 1 - F_{\Delta}(t),
\qquad
t \geq 0,
\label{A1a}
\end{equation}
where $F_{\Delta}(t)$ is given in \eqref{P3}. If we, then, differentiate \eqref{A1a} with respect to $t$, we obtain
\begin{equation}
\mathbb{E}\left[\frac{t^{D_{m,n}-1}}{\left(D_{m,n}-1\right)!} \right] e^{-t} = F_{\Delta}'(t) = : f_{\Delta}(t),
\qquad
t \geq 0.
\label{A1b}
\end{equation}

From \eqref{A2} we also have
\begin{equation}
\mathbb{E}\left[g(\Delta_{m,n}) \, | \, D_{m,n}\right] = \frac{1}{\left(D_{m,n}-1\right)!} \int_0^{\infty} g(t) \, t^{D_{m,n}-1} e^{-t} dt,
\label{A3}
\end{equation}
where $g(t)$ is any function for which the integral makes sense. Taking expectations in \eqref{A3} yields
\begin{equation}
\mathbb{E}\left[g(\Delta_{m,n})\right] = \mathbb{E}\left[\frac{1}{\left(D_{m,n}-1\right)!} \int_0^{\infty} g(t) \, t^{D_{m,n}-1} e^{-t} dt\right],
\label{A4}
\end{equation}
and if
\begin{equation}
\mathbb{E}\left[\frac{1}{\left(D_{m,n}-1\right)!} \int_0^{\infty} \left|g(t)\right| t^{D_{m,n}-1} e^{-t} dt\right] < \infty,
\label{A5}
\end{equation}
then Fubini's theorem allows us to interchange expectation and integral in \eqref{A4} and obtain
\begin{equation}
\mathbb{E}\left[g(\Delta_{m,n})\right] = \int_0^{\infty} g(t) \, \mathbb{E}\left[\frac{t^{D_{m,n}-1}}{\left(D_{m,n}-1\right)!} \right] e^{-t} dt.
\label{A6}
\end{equation}
Of course, \eqref{A6} is, also, an immediate consequence of \eqref{A1b}.

Let us look at some examples. If $g(t) = t^z$ for some complex $z$, then \eqref{A3} becomes
\begin{equation}
\mathbb{E}\left[\Delta_{m,n}^z \, | \, D_{m,n}\right] = \frac{1}{\left(D_{m,n}-1\right)!} \int_0^{\infty} t^{D_{m,n}+z-1} e^{-t} dt
= \frac{\Gamma\left(D_{m,n}+z\right)}{\left(D_{m,n}-1\right)!},
\label{A7}
\end{equation}
where $\Gamma(\cdot)$ is the Gamma function. And since $D_{m,n} \geq mn$, the above integral converges for $\Re(z) > -mn$.

In the case where $z=r$, a positive integer, formula \eqref{A7} becomes
\begin{equation*}
\mathbb{E}\left[\Delta_{m,n}^r \, | \, D_{m,n}\right] = \frac{\left(D_{m,n}+r-1\right)!}{\left(D_{m,n}-1\right)!}
= D_{m,n} \left(D_{m,n}+1\right) \cdots \left(D_{m,n} + r - 1\right),
\end{equation*}
or
\begin{equation}
\mathbb{E}\left[\Delta_{m,n}^r \, | \, D_{m,n}\right] = D_{m,n}^{(r)},
\label{A9}
\end{equation}
where we have used the notation
\begin{equation}
M^{(r)} := M (M+1) \cdots (M + r - 1),
\qquad
r \geq 1.
\label{P5}
\end{equation}
If $z = -r$, where $r = 1, 2, \ldots, (mn-1)$, formula \eqref{A7} yields
\begin{align}
\mathbb{E}\left[\Delta_{m,n}^{-r} \, | \, D_{m,n}\right] &= \frac{\left(D_{m,n}-r-1\right)!}{\left(D_{m,n}-1\right)!}
\nonumber
\\
&= \frac{1}{(D_{m,n} - r) (D_{m,n} - r + 1 ) \cdots (D_{m,n} - 1 )}.
\label{A9n}
\end{align}

Now, taking expectations in \eqref{A7} and invoking \eqref{P3} yields (after integrating by parts once)
\begin{align}
\mathbb{E}\left[\frac{\Gamma\left(D_{m,n}+z\right)}{\left(D_{m,n}-1\right)!}\right] &= \mathbb{E}\left[\Delta_{m,n}^z\right]
\nonumber
\\
&= z\int_0^{\infty} \left\{1 - \left[1 - S_{m}(t/n) e^{-t/n} \right]^n \right\} t^{z-1} dt,
\quad
\Re(z) > 0,
\nonumber
\\
&= z n^z \int_0^{\infty} \left\{1 - \left[1 - S_{m}(\tau) e^{-\tau} \right]^n \right\} \tau^{z-1} d\tau,
\quad\ \
\Re(z) > 0.
\label{A10}
\end{align}
In the case where $z = r$ is a positive integer formula \eqref{A10} becomes
\begin{equation}
\mathbb{E}\left[D_{m,n}^{(r)} \right] = \mathbb{E}\left[ \Delta_{m,n}^r \right]
= r n^r \int_0^{\infty} \left\{1 - \left[1 - S_{m}(\tau) e^{-\tau} \right]^n \right\} \tau^{r-1} d\tau.
\label{P7}
\end{equation}
The quantity $\mathbb{E}\left[ D_{m,n}^{(r)} \right]$ is called the
$r$\textit{-th rising moment of} $D_{m,n}$. For $r=1$ and $r=2$ formula \eqref{P7} gives
\begin{equation}
\mathbb{E}[D_{m,n}] = \mathbb{E}\left[ \Delta_{m,n} \right] = n \int_0^{\infty} \left\{1 - \left[1 - S_{m}(\tau) e^{-\tau} \right]^n\right\} d\tau,
\label{5}
\end{equation}
\begin{align}
\mathbb{E}\left[D_{m,n}^{(2)}\right] &= \mathbb{E}\Big[D_{m,n}\big(D_{m,n} + 1\big)\Big]
\nonumber
\\
&= \mathbb{E}\left[ \Delta_{m,n}^2 \right]
= 2 n^2 \int_0^{\infty} \left\{1 - \left[1 - S_{m}(\tau) e^{-\tau} \right]^n\right\} \tau \, d\tau
\label{5A}
\end{align}
(thus $\mathbb{V}\left[D_{m,n}\right] = \mathbb{V}\left[\Delta_{m,n}\right] - \mathbb{E}\left[\Delta_{m,n}\right]$).
Formulas \eqref{5} and \eqref{5A} were first derived in \cite{B} by a different approach.

Let us also notice that by taking expectations in \eqref{A9n} we obtain
\begin{equation}
\mathbb{E}\left[\frac{1}{(D_{m,n} - r) (D_{m,n} - r + 1 ) \cdots (D_{m,n} - 1 )}\right]
= \mathbb{E}\left[\Delta_{m,n}^{-r}\right]
\label{A9nE}
\end{equation}
for $r = 1, 2, \ldots, (mn-1)$. In particular,
\begin{equation}
\mathbb{E}\left[\frac{1}{\Delta_{m,n}}\right] = \mathbb{E}\left[\frac{1}{D_{m,n} - 1}\right]
\label{A9nEa}
\end{equation}
and
\begin{equation}
\mathbb{E}\left[\frac{1}{\Delta_{m,n}^2}\right] = \mathbb{E}\left[\frac{1}{(D_{m,n} - 2)(D_{m,n} - 1)}\right]
= \mathbb{E}\left[\frac{1}{D_{m,n} - 2}\right] - \mathbb{E}\left[\frac{1}{D_{m,n} - 1}\right].
\label{A9nEb}
\end{equation}

Suppose now that $g(t) = e^{z t}$. Then \eqref{A3} becomes
\begin{align}
\mathbb{E}\left[e^{z \Delta_{m,n}} \, | \, D_{m,n}\right] &= \frac{1}{\left(D_{m,n}-1\right)!} \int_0^{\infty} t^{D_{m,n}-1} e^{-(1-z)t} dt
\nonumber
\\
&= \frac{1}{(1-z)^{D_{m,n}}},
\qquad \Re(z) < 1
\label{A11}
\end{align}
(first we show that the second equation in \eqref{A11} holds for real $z < 1$ and then we use analytic continuation).

Hence, by taking expectations in \eqref{A11} and invoking \eqref{P3} we obtain (after integrating by parts once)
\begin{align}
\mathbb{E}\left[\frac{1}{(1-z)^{D_{m,n}}}\right] &= \mathbb{E}\left[e^{z \Delta_{m,n}}\right]
\nonumber
\\
\qquad &= 1 + z \int_0^{\infty} \left\{1 - \left[1 - S_m(t/n) e^{-t/n} \right]^n\right\} e^{zt} dt,
\quad
\Re(z) < \frac{1}{n}
\nonumber
\\
\qquad &= 1 + zn \int_0^{\infty} \left\{1 - \left[1 - S_m(\tau) e^{-\tau} \right]^n\right\} e^{nz\tau} d\tau,
\quad\
\Re(z) < \frac{1}{n}.
\label{A12}
\end{align}

\section{Asymptotics of the moments of $D_{m,n}$ as $m, n \to \infty$}
For typographical convenience we set
\begin{equation}
F_m(x) := 1 - S_m(x) e^{-x},
\qquad
x \geq 0,
\label{M0}
\end{equation}
where $S_m(\cdot)$ is given by \eqref{7}. Observe that
\begin{equation}
F_m(x) = \mathbb{P}\{\Theta_m \leq x\},
\label{M0a}
\end{equation}
where $\Theta_m$ is an Erlang random variable with parameters $m$ and $1$.

In view of \eqref{M0}, formula \eqref{A10} can be written as
\begin{equation}
\mathbb{E}\left[\frac{\Gamma\left(D_{m,n}+z\right)}{\left(D_{m,n}-1\right)!}\right] = \mathbb{E}\left[\Delta_{m,n}^z\right]
= z n^z \int_0^{\infty} \left[1 - F_m(\tau)^n \right] \tau^{z-1} d\tau,
\quad
\Re(z) > 0.
\label{M0b}
\end{equation}
In order to determine the asymptotic behavior of the integral in \eqref{M0b} one has to locate a relatively narrow interval, say $[x_1, x_2]$ of
values of $x$ in which the values of $S_m(x) e^{-x} = 1 - F_m(x)$ change from $\gg 1/n$ to $\ll 1/n$, so that $F_m(x_1)^n \to 0$ and
$F_m(x_2)^n \to 1$ as $n \to \infty$.

One observation, somehow relevant to the above task, is that (as long as $m \geq 2$) the unique point of inflection of $1 - F_m(x)$ (and of
$F_m(x)$) is located at $x = m-1$.

Now, since $\Theta_m$ of \eqref{M0a} can be expressed as a sum of $m$ independent exponential random variables with
parameter $1$, we can apply the Berry-Esseen Theorem and obtain the estimate
\begin{equation}
\left|\big[1 - F_m(x)\big] - \Phi\left(\frac{m-x}{\sqrt{m}}\right)\right| \leq \frac{C}{\sqrt{m}},
\qquad
x \in \mathbb{R}, \ \; m \in \mathbb{N} := \{1, 2, \ldots\},
\label{M1}
\end{equation}
where $\Phi(\cdot)$ is the standard normal distribution function and $C > 0$ is independent of both $x$ and $m$. For example, one implication of \eqref{M1} is
\begin{equation}
1 - F_m\big(m + O(\sqrt{m}\,)\big) \gg \frac{1}{n},
\qquad
n \to \infty.
\label{M2}
\end{equation}

In the case where $|m-x| /\sqrt{m} \to \infty$ the Berry-Esseen estimate \eqref{M1} becomes too crude. For instance, if
$(m-x) /\sqrt{m} \to -\infty$, then $\Phi\big((m-x)/\sqrt{m}\big)$ may become much smaller than the error bound $C /\sqrt{m}$.
Fortunately, in this case one can use a nice asymptotic formula due to F.G. Tricomi (1950) \cite{T} regarding the incomplete Gamma function.

First, let us notice that formula \eqref{7} implies that $S_m'(y) = S_{m-1}(y)$, hence (as we have already noticed in \eqref{A1b})
\begin{equation}
\frac{d}{dy} \left[S_m(y) e^{-y}\right] = S_{m-1}(y) e^{-y} - S_m(y) e^{-y} = -\frac{y^{m-1}}{(m-1)!} \, e^{-y}.
\label{D3}
\end{equation}
Thus, by integrating \eqref{D3} from $x$ to $\infty$ we obtain
\begin{equation}
1 - F_m(x) = S_m(x) e^{-x}
= \frac{1}{(m-1)!} \int_{x}^{\infty} y^{m-1} e^{-y} dy = \frac{\Gamma(m, x)}{(m-1)!},
\label{D4}
\end{equation}
where $\Gamma(\cdot \,, \cdot)$ is the upper incomplete Gamma function.

Tricomi \cite{T} (see, also, \cite{G}) has shown that if
\begin{equation}
x - \mu \to \infty
\qquad \text{and} \qquad
\frac{\sqrt{\mu}}{x - \mu} \to 0,
\label{D5}
\end{equation}
then
\begin{equation}
\Gamma(\mu+1, x) =
\frac{x^{\mu+1} e^{-x}}{x - \mu}
\left[1 - \frac{\mu}{(x - \mu)^2} + \frac{2\mu}{(x - \mu)^3} + O\left(\frac{\mu^2}{(x - \mu)^4}\right)\right].
\label{D6}
\end{equation}
Formula \eqref{D6} holds for complex $\mu$ and $x$ as long as the argument of the quantity $\sqrt{\mu}/(x - \mu)$ ultimately remains between
$-3\pi/4$ and $3\pi/4$. In fact, Tricomi \cite{T} has given the complete asymptotic expansion of $\Gamma(\mu+1, x)$, but, for our
purposes formula \eqref{D6} is more than sufficient. Indeed, by using \eqref{D6} in \eqref{D4} with $\mu = m-1$ and invoking Stirling's formula we can conclude that
\begin{align}
1 - F_m(x) &= S_m(x) e^{-x}
\nonumber
\\
&= \frac{1}{\sqrt{2\pi}} \, e^{-(x-m)} \left(\frac{x}{m}\right)^m
\frac{\sqrt{m}}{x - m + 1}
\left[1 + O\left(\frac{m}{(x - m)^2}\vee\frac{1}{m}\right)\right],
\label{D6a}
\end{align}
provided that \eqref{D5} is satisfied, namely
\begin{equation}
x - m \to \infty
\qquad \text{and} \qquad
\frac{\sqrt{m}}{x - m} \to 0
\label{D5a}
\end{equation}
(in \eqref{D6a} we have used the notation $a \vee b$ for the maximum of $a$ and $b$).

\subsection{The ``supercritical" case $m \gg \ln n$}
In the case where $m$ grows faster than $\ln n$, it turns out that the desired jump of $F_m(x)^n$ from $0$ to $1$ happens in the interval
$m \leq x \leq m[1 + \varepsilon(n)]$, where $\varepsilon(n)$ is given by \eqref{M22} below. This is the main ingredient in the proof of the theorem that follows.

\medskip

\textbf{Theorem 1.} Let $z$ be a fixed complex number with $\Re(z) > 0$ and assume that $m = m(n)$ is such that $m(n) \gg \ln n$, i.e.
\begin{equation}
\frac{\ln n}{m(n)} \to 0
\qquad \text{as }\;
n \to \infty.
\label{M4a}
\end{equation}
Then
\begin{equation}
\mathbb{E}\left[\frac{\Gamma\left(D_{m,n}+z\right)}{\left(D_{m,n}-1\right)!}\right] = \mathbb{E}\left[\Delta_{m,n}^z\right]
\sim n^z m^z,
\qquad
n \to \infty
\label{M4}
\end{equation}
(recall that the variable $\Delta_{m,n}$ has been introduced in Subsection 1.1).

\smallskip

\textit{Proof}. The substitution $\tau = m \xi$ in the integral of \eqref{M0b} yields
\begin{equation}
\mathbb{E}\left[\frac{\Gamma\left(D_{m,n}+z\right)}{\left(D_{m,n}-1\right)!}\right] = \mathbb{E}\left[\Delta_{m,n}^z\right] = z n^z m^z I(z),
\label{M5a}
\end{equation}
where
\begin{equation}
I(z) := \int_0^{\infty} \left[1 - F_m(m \xi)^n \right] \xi^{z-1} d\xi.
\label{M5b}
\end{equation}
We split $I(z)$ as
\begin{equation}
I(z) = I_0(z) + I_1(z),
\label{M6a}
\end{equation}
where
\begin{equation}
I_0(z) := \int_0^{1+\varepsilon} \left[1 - F_m(m \xi)^n \right] \xi^{z-1} d\xi,
\quad
I_1(z) := \int_{1+\varepsilon}^{\infty} \left[1 - F_m(m \xi)^n \right] \xi^{z-1} d\xi,
\label{M6b}
\end{equation}
where $\varepsilon = \varepsilon(n)$ is positive and approaches $0$ as $n \to \infty$. The specific form of $\varepsilon(n)$ will be decided later in the proof.

Using the Berry-Esseen estimate \eqref{M1} and applying bounded convergence to the first integral in \eqref{M6b} we obtain immediately that
\begin{equation}
\lim_n I_0(z) = \int_0^1 \xi^{z-1} d\xi = \frac{1}{z}.
\label{M3}
\end{equation}
Therefore, in view of \eqref{M5a}, \eqref{M5b}, \eqref{M6a}, \eqref{M6b}, and \eqref{M3}, the proof of \eqref{M4} will be completed if we show
that we can choose $\varepsilon = \varepsilon(n)$ so that
\begin{equation}
\lim_n I_1(z) = 0.
\label{M3a}
\end{equation}

From \eqref{M6b} we have
\begin{equation}
\left|I_1(z)\right| \leq \int_{1+\varepsilon}^{\infty} \left[1 - F_m(m \xi)^n \right] \xi^{\sigma-1} d\xi,
\qquad \text{where }\;
\sigma := \Re(z) > 0.
\label{M7}
\end{equation}
Now, in view of \eqref{M0},
\begin{equation}
F_m(m \xi)^n = \left[1 - S_m(m \xi) e^{-m \xi}\right]^n \geq 1 - n  S_m(m \xi) e^{-m \xi}.
\label{M8}
\end{equation}
Hence, \eqref{M7} implies
\begin{equation}
\left|I_1(z)\right| \leq n\int_{1+\varepsilon}^{\infty} S_m(m \xi) e^{-m \xi} \xi^{\sigma-1} d\xi.
\label{M9}
\end{equation}
Next we notice that for $\xi \geq 1$ formula \eqref{7} implies easily that
\begin{equation}
S_m(m \xi) \leq m \frac{(m \xi)^{m-1}}{(m-1)!} = \frac{m^{m+1}}{m!} \, \xi^{m-1}.
\label{M10}
\end{equation}
Thus, by using \eqref{M10} in \eqref{M9} we get
\begin{equation}
\left|I_1(z)\right| \leq n \frac{m^{m+1}}{m!} \int_{1+\varepsilon}^{\infty} \xi^{\sigma-2} \xi^m e^{-m \xi} d\xi
= n \frac{m^{m+1}}{m!} \, e^{-m} \Lambda,
\label{M11}
\end{equation}
where we have set
\begin{equation}
\Lambda := \int_{\varepsilon}^{\infty} (1+t)^{\sigma-2} e^{-m[t - \ln(1+t)]} dt.
\label{M12}
\end{equation}
For the function
\begin{equation}
\rho(t) := t - \ln(1+t)
\label{M13}
\end{equation}
which appears in the exponent of the integrand in \eqref{M12} we have
\begin{equation}
\rho(0) = \rho'(0) = 0
\qquad \text{and} \qquad
\rho'(t) = 1 - \frac{1}{t+1} > 0 \ \; \text{for } t \in (0, \infty),
\label{M14}
\end{equation}
so that $\rho(t)$ is strictly increasing on $[0, \infty)$. Therefore, for the integral in \eqref{M12} we have the estimate
\begin{equation}
\Lambda < 2 \int_{\varepsilon}^1 (1+t)^{\sigma-2} e^{-m[t - \ln(1+t)]} dt < 2^{\sigma-1} \int_{\varepsilon}^1 e^{-m[t - \ln(1+t)]} dt,
\label{M15}
\end{equation}
as long as $m$ is sufficiently large.

Now, from \eqref{M14} (and the fact that $\rho''(0) = 1$) it also follows that for all sufficiently small $\delta > 0$ we must have
\begin{equation}
\rho(t) \geq \delta t^2
\qquad \text{for every }\; t \in [0, 1].
\label{M16}
\end{equation}
If we choose such a $\delta$, then \eqref{M15} implies
\begin{equation}
\Lambda < 2^{\sigma-1} \int_{\varepsilon}^1 e^{-m \delta t^2} dt < 2^{\sigma-1} \int_{\varepsilon}^{\infty} e^{-m \delta t^2} dt
= \frac{2^{\sigma-1}}{\sqrt{\delta m}} \int_{\varepsilon\sqrt{\delta m}}^{\infty} e^{-x^2} dx.
\label{M17}
\end{equation}
Thus, in view of \eqref{M11} and \eqref{M17} (and Stirling's formula) the limit \eqref{M3a} will hold if we can find an
$\varepsilon = \varepsilon(n)$ (with $\varepsilon(n) \to 0$ as $n \to \infty$) such that
\begin{equation}
\int_{\varepsilon\sqrt{\delta m}}^{\infty} e^{-x^2} dx = o\left(\frac{1}{n}\right),
\qquad
n \to \infty.
\label{M18}
\end{equation}
In order to satisfy \eqref{M18} it is clearly necessary to have
\begin{equation}
\varepsilon\sqrt{m} \to \infty
\qquad \text{as }\;
n \to \infty.
\label{M19}
\end{equation}
Now, if \eqref{M19} is satisfied, then it is not hard to see (e.g., by applying L'H\^{o}pital's rule) that
\begin{equation}
\int_{\varepsilon\sqrt{\delta m}}^{\infty} e^{-x^2} dx \sim \frac{1}{2\varepsilon\sqrt{\delta m}}\,e^{-m \delta\varepsilon^2},
\label{M20}
\end{equation}
Hence, under \eqref{M19}, formula \eqref{M18} is equivalent to
\begin{equation}
\frac{1}{\varepsilon\sqrt{m}}\,e^{-m \delta\varepsilon^2} = o\left(\frac{1}{n}\right),
\qquad
n \to \infty.
\label{M21}
\end{equation}
Notice that if $m = m(n)$ does not grow faster than $\ln n$, then it is not possible to simultaneously satisfy $\varepsilon(n) \to 0$, \eqref{M19},
and \eqref{M21}. This indicates that without the assumption \eqref{M4a} formula \eqref{M4} may not hold.

To satisfy \eqref{M21} we can pick a $\kappa > 1$ and then take
\begin{equation}
\varepsilon = \varepsilon(n) = \sqrt{\frac{\kappa \ln n}{\delta m}}.
\label{M22}
\end{equation}
Since \eqref{M18} is equivalent to \eqref{M21}, we can conclude that \eqref{M18} is satisfied if $\varepsilon(n)$ is chosen as above
(also it is clear that, under assumption \eqref{M4a}, the above choice of $\varepsilon(n)$ satisfies $\lim_n \varepsilon(n) = 0$ as well as \eqref{M19}).
\hfill $\blacksquare$

\medskip

In the case where $z = r \in \mathbb{N}$, formula \eqref{M4} becomes
\begin{equation}
\mathbb{E}\left[D_{m,n}^{(r)}\right] = \mathbb{E}\left[\Delta_{m,n}^r\right]
\sim n^r m^r,
\qquad
n \to \infty,
\label{M23}
\end{equation}
from which it follows that
\begin{equation}
\mathbb{E}\left[D_{m,n}^r\right] \sim n^r m^r,
\qquad
n \to \infty.
\label{M24}
\end{equation}
In particular,
\begin{equation}
\mathbb{E}\left[D_{m,n}\right] \sim n m,
\qquad
n \to \infty,
\label{M25}
\end{equation}
which is in agreement with the corresponding result which appeared in \cite{Sh} and \cite{S-H}.

Notice that, roughly speaking, formula \eqref{M25} tells us that, on the average, if all $n$ coupons have already been detected $m$ times, where
$m \gg \ln n$, then each additional detection (of all coupons) ``costs" $n$. Apart from being interesting by itself, formula \eqref{M25} is used in
the proof of Theorem 4 below.

\medskip

\textbf{Remark 1.} Suppose $n$ is fixed. Then formula \eqref{M4} is still valid, where now $m \to \infty$, i.e. if $\Re(z) > 0$, then
\begin{equation}
\mathbb{E}\left[\frac{\Gamma\left(D_{m,n}+z\right)}{\left(D_{m,n}-1\right)!}\right] = \mathbb{E}\left[\Delta_{m,n}^z\right]
\sim n^z m^z,
\qquad
m \to \infty;
\label{R1M4}
\end{equation}
in particular, for $r \in \mathbb{N}$ we get
\begin{equation}
\mathbb{E}\left[D_{m,n}^{(r)}\right] = \mathbb{E}\left[\Delta_{m,n}^r\right]
\sim n^r m^r,
\qquad
m \to \infty.
\label{R1M23}
\end{equation}

We can prove \eqref{R1M4} by slightly modifying the proof of Theorem 1. Here is the key element in the modified proof: Instead of \eqref{M22} we
now pick a (constant, but otherwise arbitrary) $\varepsilon > 0$. The Berry-Esseen
theorem implies (as $m \to \infty$)
\begin{equation}
F_m(m \xi) = \Phi\left(\frac{m \xi - m}{\sqrt{m}}\right) + O\left(\frac{1}{\sqrt{m}}\right) = \Phi\left((\xi - 1)\sqrt{m}\right)
+ O\left(\frac{1}{\sqrt{m}}\right),
\label{R1M1}
\end{equation}
uniformly in $\xi$. Thus,
\begin{equation}
F_m(m \xi)^n = O\left(\frac{1}{\sqrt{m}}\right),
\qquad\text{if }\;
0 \leq \xi \leq 1 - \varepsilon ,
\label{R1M2a}
\end{equation}
while
\begin{equation}
F_m(m \xi)^n = 1 + O\left(\frac{1}{\sqrt{m}}\right),
\qquad\text{if }\;
\xi \geq 1 + \varepsilon.
\label{R1M2b}
\end{equation}

\subsection{The critical case $m \sim \beta\ln n$}
In the case where $m$ grows like $\ln n$, the formula \eqref{M4} is no longer true. As we will see in the following lemma, the reason is that,
when $m \sim \beta\ln n$ for some $\beta > 0$,
in contrast with the supercritical case, the values of $x$ at which $S_m(x) e^{-x} = 1 - F_m(x)$ changes from $\gg 1/n$ to $\ll 1/n$ are quite away
from $m$.

\medskip

\textbf{Lemma 1.} Suppose $m = m(n)$, $n = 1, 2, \ldots$, is a given sequence of positive integers such that
\begin{equation}
m(n) \sim \beta \ln n,
\qquad
n \to \infty,
\label{CL1}
\end{equation}
where $\beta > 0$ is a fixed constant, and let $\alpha$ be the unique solution of the equation
\begin{equation}
\alpha - \beta\ln \alpha = \beta - \beta \ln\beta + 1
\label{CL2}
\end{equation}
in the interval $(\beta, \infty)$ (thus $\alpha > \beta$). Then, there are two sequences $a_1 = a_1(n) \to 0$ and $a_2 = a_2(n) \to 0$, with
$a_1(n) < a_2(n)$ such that (as $n \to \infty$):
\begin{equation}
\text{If}\quad
\chi_1 = \chi_1(n) := \alpha \ln n + a_1(n) \ln n,
\qquad \text{then}\quad
S_m(\chi_1) e^{-\chi_1} \gg \frac{1}{n},
\label{CL3a}
\end{equation}
while
\begin{equation}
\text{if}\quad
\chi_2 = \chi_2(n) := \alpha \ln n + a_2(n) \ln n,
\qquad \text{then}\quad
S_m(\chi_2) e^{-\chi_2} \ll \frac{1}{n}.
\label{CL3b}
\end{equation}

\smallskip

\textit{Proof}. Notice that, in view of \eqref{CL1}, \eqref{CL2}, \eqref{CL3a}, and \eqref{CL3b}, if $x = \chi_1$ or $x = \chi_2$, then
the conditions \eqref{D5a} are satisfied. Thus, \eqref{D6a} implies
\begin{equation}
S_m(\chi_j) e^{-\chi_j} = \frac{1}{\sqrt{2\pi}} \frac{e^{-(\chi_j -m)}}{\sqrt{m}} \left(\frac{\chi_j}{m}\right)^m \frac{m}{\chi_j - m}
\left[1 + O\left(\frac{1}{\ln n}\right)\right],
\qquad
j = 1, 2.
\label{CL4}
\end{equation}
We will now analyze each of the three factors appearing in the right-hand side of \eqref{CL4}. But, before we start, let us express assumption
\eqref{CL1} in the equivalent form
\begin{equation}
m(n) := \beta \ln n + b(n) \ln n,
\qquad\text{where }\; b(n) \to 0
\label{CL1a}
\end{equation}
(the sequence $b = b(n)$ is assumed given).

For the factor $m / (\chi_j - m)$ it follows immediately that
\begin{equation}
\frac{m}{\chi_j - m} \sim \frac{\beta}{\alpha - \beta}.
\label{CL5}
\end{equation}

For the factor $e^{-(\chi_j -m)}/\sqrt{m}$ we have, in view of \eqref{CL3a}, \eqref{CL3b}, and \eqref{CL1a},
\begin{equation}
\frac{e^{-(\chi_j -m)}}{\sqrt{m}} \sim \frac{1}{\sqrt{\beta}} \, \frac{1}{\sqrt{\ln n}} \, \frac{1}{n^{\alpha-\beta}} \, e^{-(a_j - b) \ln n}.
\label{CL6}
\end{equation}

Finally, we analyze the factor $(\chi_j/m)^m$:
\begin{align}
\left(\frac{\chi_j}{m}\right)^m &= \left(\frac{\alpha + a_j}{\beta + b}\right)^m
= \left(\frac{\alpha}{\beta}\right)^m \left(\frac{1 + a_j/\alpha}{1+ b/\beta}\right)^m
\nonumber
\\
&= \left(\frac{\alpha}{\beta}\right)^{\beta \ln n + b \ln n} \left(\frac{1 + a_j/\alpha}{1+ b/\beta}\right)^{\beta \ln n + b \ln n}
\nonumber
\\
&= \frac{e^{\,(\ln\alpha - \ln\beta)\,b\ln n}}{n^{\beta\ln\beta - \beta\ln\alpha}} \,
e^{\,\left[\ln(1 + a_j/\alpha) - \ln(1+ b/\beta)\right](\beta \ln n + b \ln n)}.
\label{CL7}
\end{align}
If we now use \eqref{CL5}, \eqref{CL6}, and \eqref{CL7} in \eqref{CL4} and invoke \eqref{CL2}, we obtain
\begin{equation}
S_m(\chi_j) e^{-\chi_j} \sim \frac{1}{\sqrt{2\pi}} \, \frac{\sqrt{\beta}}{\alpha - \beta} \, \frac{e^{A_j(n) \ln n}}{n},
\qquad
j = 1, 2,
\label{CL8}
\end{equation}
where for typographical convenience we have set
\begin{align}
A_j(n) &:=(\ln\alpha - \ln\beta + 1)\,b - (\beta + b) \ln\left(1+ \frac{b}{\beta}\right)
\nonumber
\\
&- a_j + (\beta + b) \ln\left(1 + \frac{a_j}{\alpha}\right) - \frac{1}{2} \frac{\ln\ln n}{\ln n},
\qquad
j = 1, 2.
\label{CL8a}
\end{align}
Our assumption $b = b(n) \to 0$, as $n \to \infty$, implies immediately that
\begin{equation}
B = B(n) := (\ln\alpha - \ln\beta + 1)\,b - (\beta + b) \ln\left(1 + \frac{b}{\beta}\right) \to 0.
\label{CL9}
\end{equation}
For the sequences $a_1$ and $a_2$ we need $a_j \to 0$. Thus,
\begin{equation}
- a_j + (\beta + b) \ln\left(1 + \frac{a_j}{\alpha}\right) = -\left(1 - \frac{\beta}{\alpha}\right) a_j + o\left(a_j\right),
\qquad
n \to \infty
\label{CL10}
\end{equation}
(recall that $0 < \beta/\alpha < 1$).

Substituting \eqref{CL9} and \eqref{CL10} in \eqref{CL8a} yields
\begin{equation}
A_j(n) = B(n) - \frac{1}{2} \frac{\ln\ln n}{\ln n} -\left(1 - \frac{\beta}{\alpha}\right) a_j + o\left(a_j\right)
\label{CL11}
\end{equation}
and from formula \eqref{CL11} it is clear that, since $B(n)$ is a given sequence with $B(n) \to 0$, we can choose a sequence $a_1 = a_1(n)$,
with $a_1(n) \to 0$ so that $A_1(n) \ln n \to \infty$. For example, just take
\begin{equation*}
a_1(n) = -\frac{2\alpha}{\alpha-\beta}\left(|B(n)| + \frac{\ln\ln n}{2\ln n}\right).
\end{equation*}
Likewise, we can choose a sequence $a_2 = a_2(n)$,
with $a_2(n) \to 0$ so that $A_2(n) \ln n \to -\infty$. For example,
\begin{equation*}
a_2(n) = \frac{2\alpha}{\alpha-\beta}\left(|B(n)| + \frac{\ln\ln n}{2\ln n}\right).
\end{equation*}
Therefore, in view of \eqref{CL8}, we have demonstrated that there are sequences $a_1(n) \to 0$ and $a_2(n) \to 0$ which satisfy \eqref{CL3a} and
\eqref{CL3b} respectively.
\hfill $\blacksquare$

\medskip

\textbf{Remark 2.} Formula \eqref{CL2} can be written as $\alpha = \beta + 1 + \beta\ln(\alpha/\beta)$ and from this, together with the fact that
$\alpha > \beta > 0$, it is obvious that $\alpha > \beta + 1$. Also, $\alpha$ increases with $\beta$ (as it is easy to check that
$d\alpha/d\beta > 0$) and $\alpha \to 1^+$ as $\beta \to 0^+$. Thus, $\alpha$ can take any value in $(1, \infty)$.

\medskip

It is remarkable that the equation \eqref{CL2} also appears, in a different context, in Sharif and Hassibi \cite{S-H2}.

\medskip

\textbf{Remark 3.} It is rather obvious that if $\tilde{a}_1(n)$ is a sequence such that $\tilde{a}_1(n) \to 0$ and $\tilde{a}_1(n) \leq a_1(n)$,
where $\chi_1(n) := \alpha \ln n + a_1(n) \ln n$ satisfies \eqref{CL3a}, then $\tilde{\chi}_1(n) := \alpha \ln n + \tilde{a}_1(n) \ln n$ also
satisfies \eqref{CL3a}. Likewise, if $\tilde{a}_2(n)$ is a sequence such that $\tilde{a}_2(n) \to 0$ and $\tilde{a}_2(n) \geq a_2(n)$,
where $\chi_2(n) := \alpha \ln n + a_2(n) \ln n$ satisfies \eqref{CL3b}, then $\tilde{\chi}_2(n) := \alpha \ln n + \tilde{a}_2(n) \ln n$ also
satisfies \eqref{CL3b}.

\medskip

We are now ready to give the (leading) asymptotic behavior of the moments of $D_{m,n}$ as $n \to \infty$ in the critical case.

\medskip

\textbf{Theorem 2.} Let $z$ be a fixed complex number with $\Re(z) > 0$ and assume that
\begin{equation}
m = m(n) \sim \beta\ln n
\qquad\text{for some constant }\;
\beta > 0.
\label{T2M0}
\end{equation}
Then,
\begin{equation}
\mathbb{E}\left[\frac{\Gamma\left(D_{m,n}+z\right)}{\left(D_{m,n}-1\right)!}\right] = \mathbb{E}\left[\Delta_{m,n}^z\right]
\sim \alpha^z n^z (\ln n)^z \sim \left(\frac{\alpha}{\beta}\right)^z n^z m^z,
\qquad
n \to \infty,
\label{T2M1}
\end{equation}
where $\alpha$ is the unique solution of the equation \eqref{CL2} in the interval $(\beta, \infty)$.

\smallskip

\textit{Proof}. We will prove the theorem by adapting the proof of Theorem 1.

From formula \eqref{M0b} we have
\begin{equation}
\mathbb{E}\left[\frac{\Gamma\left(D_{m,n}+z\right)}{\left(D_{m,n}-1\right)!}\right] = \mathbb{E}\left[\Delta_{m,n}^z\right] = n^z m^z J(z),
\label{T2M2a}
\end{equation}
where
\begin{equation}
J(z) := \int_0^{\infty} \left[1 - F_m(m \xi)^n \right] z\xi^{z-1} d\xi.
\label{T2M2b}
\end{equation}
This time we split $J(z)$ as
\begin{equation}
J(z) = J_1(z) + J_2(z) + J_3(z),
\label{T2M3}
\end{equation}
where
\begin{equation}
J_1(z) := \int_0^{\chi_1/m} \left[1 - F_m(m \xi)^n \right] z\xi^{z-1} d\xi,
\label{T2M3a}
\end{equation}
\begin{equation}
J_2(z) := \int_{\chi_1/m}^{\chi_2/m} \left[1 - F_m(m \xi)^n \right] z\xi^{z-1} d\xi,
\label{T2M3b}
\end{equation}
and
\begin{equation}
J_3(z) := \int_{\chi_2/m}^{\infty} \left[1 - F_m(m \xi)^n \right] z\xi^{z-1} d\xi,
\label{T2M3c}
\end{equation}
where $\chi_1$ and $\chi_2$ are as in \eqref{CL3a} and \eqref{CL3b} respectively. Notice that \eqref{CL3a}, \eqref{CL3b},
and \eqref{T2M0} imply that
\begin{equation}
\frac{\chi_j}{m} = \frac{\alpha}{\beta} + o(1)
\qquad \text{as }\; n \to \infty.
\label{T2M4}
\end{equation}

From the display \eqref{CL3a} and the fact that $S_m(x)e^{-x}$ is decreasing on $[0, \infty)$ we get (in view of \eqref{T2M4}) that
\begin{equation}
n S_m(m \xi) e^{-m\xi} \to \infty,
\qquad \text{if }\;
\xi \in \left[0, \frac{\alpha}{\beta} - \varepsilon\right],
\label{T2M5}
\end{equation}
for any given $\varepsilon > 0$. Consequently,
\begin{equation}
\lim_n F_m(m \xi)^n = \lim_n \left[1 - S_m(m \xi) e^{-m\xi}\right]^n = 0
\qquad \text{for all }\;
\xi \in \left[0, \frac{\alpha}{\beta} - \varepsilon\right].
\label{T2M6}
\end{equation}
Let us break the integral $J_1(z)$ of \eqref{T2M3a} as
\begin{equation}
J_1(z) = \int_0^{\frac{\alpha}{\beta} - \varepsilon} \left[1 - F_m(m \xi)^n \right] z\xi^{z-1} d\xi
+ \int_{\frac{\alpha}{\beta} - \varepsilon}^{\frac{\chi_j}{m}} \left[1 - F_m(m \xi)^n \right] z\xi^{z-1} d\xi.
\label{T2M7}
\end{equation}
By \eqref{T2M4} and the fact that $F_m(x)$ is a distribution function, the second integral in the right-hand side of \eqref{T2M7} is bounded by
$|z| (\alpha/\beta)^{\sigma-1} [\varepsilon + o(1)]$, where, as in \eqref{M7}, $\sigma = \Re(z)$. Therefore, in view of \eqref{T2M6} and the fact that $\varepsilon$ is arbitrary, we can conclude from \eqref{T2M7} that
\begin{equation}
J_1(z) \to \int_0^{\frac{\alpha}{\beta}} z\xi^{z-1} d\xi = \left(\frac{\alpha}{\beta}\right)^z
\qquad \text{as }\; n \to \infty.
\label{T2M8}
\end{equation}
Now, the integral $J_2(z)$ of \eqref{T2M3b} is very easy to treat. We have
\begin{equation}
\left|J_2(z)\right| \leq |z|\int_{\frac{\alpha}{\beta} + o(1)}^{\frac{\alpha}{\beta} + o(1)} \xi^{\sigma-1} d\xi,
\label{T2M9}
\end{equation}
which implies immediately that
\begin{equation}
J_2(z) \to 0
\qquad \text{as }\; n \to \infty.
\label{T2M10}
\end{equation}
Thus, in view of \eqref{T2M2a}, \eqref{T2M2b}, \eqref{T2M3}, \eqref{T2M3a}, \eqref{T2M3b}, \eqref{T2M3c}, \eqref{T2M8}, and \eqref{T2M9}, in order
to complete the proof of \eqref{T2M1} it remains to show that we can choose $a_2 = a_2(n)$ (recall \eqref{CL3b} and Remark 3) so that
\begin{equation}
\lim_n J_3(z) = 0.
\label{T2M11a}
\end{equation}
In the same way we derived formula \eqref{M11} in the proof of Theorem 1, we can now get
\begin{equation}
\left|J_3(z)\right| \leq |z| n \frac{m^{m+1}}{m!} \int_{\frac{\alpha}{\beta} + a_2}^{\infty} \xi^{\sigma-2} \xi^m e^{-m \xi} d\xi,
\qquad\text{with }\;
\sigma = \Re(z).
\label{T2M11}
\end{equation}
Since $z$ is fixed, application of Stirling's formula in \eqref{T2M11} yields
\begin{equation}
\left|J_3(z)\right| \leq C n e^m \sqrt{m} \int_{\frac{\alpha}{\beta} + a_2}^{\infty} \xi^{\sigma-2} e^{-m (\xi - \ln\xi)} d\xi
\qquad\text{for some constant }\: C > 0.
\label{T2M12}
\end{equation}
Let us set
\begin{equation}
\phi(\xi) := \xi - \ln\xi.
\label{T2M13}
\end{equation}
Recall that $\alpha/\beta > 1$, while $a_2 = a_2(n)$ (see Lemma 1) is a sequence approaching $0$, which, without loss of generality (in view of
Remark 3), can be assumed positive. And since $\phi'(\xi) = 1 - \xi^{-1}$ and $\phi''(\xi) = \xi^{-2}$, it follows that $\phi(\xi)$ is convex and
increasing on $[\alpha/\beta + a_2, \infty)$. Therefore,
\begin{equation}
\phi(\xi) \geq \phi\left(\frac{\alpha}{\beta} + a_2\right)
+ \phi'\left(\frac{\alpha}{\beta} + a_2\right) \left(\xi - \frac{\alpha}{\beta} - a_2\right)
\qquad\text{for }\;
\xi \geq \frac{\alpha}{\beta} + a_2,
\label{T2M14}
\end{equation}
and, hence, by using \eqref{T2M13} and \eqref{T2M14} in \eqref{T2M12} we obtain
\begin{equation*}
\left|J_3(z)\right| \leq C n e^m \sqrt{m} \int_{\frac{\alpha}{\beta} + a_2}^{\infty} \xi^{\sigma-2}
e^{-m\left[\phi\left(\frac{\alpha}{\beta} + a_2\right)
+ \phi'\left(\frac{\alpha}{\beta} + a_2\right) \left(\xi - \frac{\alpha}{\beta} - a_2\right)\right]} d\xi,
\end{equation*}
or
\begin{align}
\left|J_3(z)\right| &\leq C n e^{m\left[1 -\phi\left(\frac{\alpha}{\beta} + a_2\right)\right]} \sqrt{m}
\int_0^{\infty} \left(t + {\frac{\alpha}{\beta} + a_2}\right)^{\sigma-2}
e^{-m \phi'\left(\frac{\alpha}{\beta} + a_2\right)t} dt
\nonumber
\\
&\leq K n e^{m\left[1 -\phi\left(\frac{\alpha}{\beta} + a_2\right)\right]} \sqrt{m}
\int_0^{\infty} e^{-m \phi'\left(\frac{\alpha}{\beta} + a_2\right)t} dt
= K \frac{n e^{m\left[1 -\phi\left(\frac{\alpha}{\beta} + a_2\right)\right]}}{\sqrt{m} \, \phi'\left(\frac{\alpha}{\beta} + a_2\right)}
\label{T2M15}
\end{align}
for some constant $K > 0$.

Now, in view of \eqref{T2M13} and the fact that $a_2 \to 0$ we have
\begin{align}
\phi\left(\frac{\alpha}{\beta} + a_2\right) &= \frac{\alpha}{\beta} + a_2 - \ln\left(\frac{\alpha}{\beta} + a_2\right)
= \frac{\alpha}{\beta} + a_2 - \ln\Bigg(\frac{\alpha}{\beta}\left(1 + \frac{\beta}{\alpha}\,a_2\right)\Bigg)
\nonumber
\\
&= \frac{\alpha}{\beta} - \ln\left(\frac{\alpha}{\beta}\right) + a_2 - \frac{\beta}{\alpha}\,a_2 + O\left(a_2^2\right)
\nonumber
\\
&= \frac{\alpha - \beta\ln\alpha + \beta\ln\beta}{\beta}  + \left(1 - \frac{\beta}{\alpha}\right)a_2 + O\left(a_2^2\right)
\nonumber
\\
&= 1 + \frac{1}{\beta}  + \left(1 - \frac{\beta}{\alpha}\right)a_2 + O\left(a_2^2\right),
\label{T2M16}
\end{align}
where the last equality follows from \eqref{CL2}.
Also,
\begin{equation}
\phi'\left(\frac{\alpha}{\beta} + a_2\right) = 1 - \left(\frac{\alpha}{\beta} + a_2\right)^{-1}
= 1 - \frac{\beta}{\alpha} \, \frac{1}{1 + \frac{\beta}{\alpha}\,a_2} = 1 - \frac{\beta}{\alpha} + O\left(a_2\right).
\label{T2M17}
\end{equation}
Substituting \eqref{T2M16} and \eqref{T2M17} in \eqref{T2M15} yields
\begin{equation}
\left|J_3(z)\right| \leq K \frac{n e^{-m\left[\frac{1}{\beta}  + \left(1 - \frac{\beta}{\alpha}\right)a_2 + O\left(a_2^2\right)\right]}}
{\sqrt{m} \left[1 - \frac{\beta}{\alpha} + O\left(a_2\right)\right]}.
\label{T2M18}
\end{equation}
Let us recall that $0 < \beta/\alpha < 1$ and $m = m(n) = \beta \ln n + b(n) \ln n$, with $b(n) \to 0$. Thus, the denominator of the fraction in
the right-hand side of \eqref{T2M18} approaches $\infty$ as $n \to \infty$. As for the numerator of that fraction, if we choose a sequence $a_2$ so
that
\begin{equation}
a_2(n) \gg |b(n)| + \frac{\ln\ln n}{\ln n}
\label{T2M19}
\end{equation}
(this choice is legitimate in view of Remark 3), then
\begin{equation}
n e^{-m\left[\frac{1}{\beta}  + \left(1 - \frac{\beta}{\alpha}\right)a_2 + O\left(a_2^2\right)\right]}
= n e^{-\left[1 + \beta\left(1 - \frac{\beta}{\alpha}\right)a_2 + o\left(a_2\right)\right] \ln n}
= e^{-\left[\beta\left(1 - \frac{\beta}{\alpha}\right)a_2 + o\left(a_2\right)\right] \ln n}
\label{T2M20}
\end{equation}
and, hence,
\begin{equation}
n e^{-m\left[\frac{1}{\beta}  + \left(1 - \frac{\beta}{\alpha}\right)a_2 + O\left(a_2^2\right)\right]} \to 0
\qquad\text{as }\; n \to \infty.
\label{T2M21}
\end{equation}
Finally, by using \eqref{T2M21} in \eqref{T2M18} we deduce that $J_3(z) \to 0$ as $n \to \infty$.
Therefore, there is a sequence $a_2(n)$ for which \eqref{T2M11a} is satisfied, and the proof is finished.
\hfill $\blacksquare$

\medskip

In the case where $z = r \in \mathbb{N}$, formula \eqref{T2M1} becomes
\begin{equation}
\mathbb{E}\left[D_{m,n}^{(r)}\right] = \mathbb{E}\left[\Delta_{m,n}^r\right]
\sim \alpha^r n^r (\ln n)^r \sim \left(\frac{\alpha}{\beta}\right)^r n^r m^r,
\qquad
n \to \infty,
\label{MM23}
\end{equation}
from which it follows that
\begin{equation}
\mathbb{E}\left[D_{m,n}^r\right] \sim \alpha^r n^r (\ln n)^r \sim \left(\frac{\alpha}{\beta}\right)^r n^r m^r,
\qquad
n \to \infty.
\label{MM24}
\end{equation}
In particular,
\begin{equation}
\mathbb{E}\left[D_{m,n}\right] \sim \alpha \, n \ln n \sim \frac{\alpha}{\beta} \, n m,
\qquad
n \to \infty,
\label{MM25}
\end{equation}
which is in agreement with the corresponding result appeared in \cite{Sh} and \cite{S-H}.
Roughly speaking, formula \eqref{MM25} tells us that, on the average, if all $n$ coupons have already been detected $m$ times, where
$m \sim \beta\ln n$, then each additional detection (of all coupons) ``costs" $(\alpha/\beta) n$. Apart from being interesting by itself, formula
\eqref{MM25} is used in the proof of Theorem 6 below.

\section{The limiting distribution of $D_{m,n}$}
Let us mention again that the results of this section are not new, since they were derived by a more direct method and in a more general setup by
G.I. Ivchenko in \cite{I1} and \cite{I2}.

Our strategy for determining the limiting distribution of $D_{m,n}$ can be described as follows.
We start with the observation that formulas \eqref{P4} and \eqref{P7} hint that, under a suitable normalization the limiting distributions of
$\Delta_{m,n}$ and $D_{m,n}$ should coincide. Hence, we can first try to find the limiting distribution of $\Delta_{m,n}$, which, thanks to
\eqref{P3}, seems an easier problem, and from that obtain the limiting distribution of $D_{m,n}$
(with the help of the ``Converging Together Lemma" --- see below).
In order, though, to determine the limiting
distribution of $\Delta_{m,n}$, we first need to come up with its correct normalization, and this, in view of \eqref{P3}, can be accomplished,
if we manage to find an expression $t$ for which $S_m(t/n) e^{-t/n} \sim Q/n$, where $Q$ is some quantity which is independent of $n$. This task
is more delicate than the one of the previous section, where we had to determine an interval of values of $x$ in which $S_m(x) e^{-x}$ changes from
$\gg 1/n$ to $\ll 1/n$. For this reason one expects that, in order to achieve the desired asymptotics, i.e. $S_m(t/n) e^{-t/n} \sim Q/n$, we may need
to impose some mild restrictions on $m(n)$.

Another thing worth repeating here is that, in order to obtain the limiting distribution of $D_{m,n}$ from the limiting distribution of
$\Delta_{m,n}$ via the Converging Together Lemma, it is necessary to have an estimate for the growth of $\mathbb{E}\left[D_{m,n}\right]$. Hence,
formulas \eqref{M25} and \eqref{MM25} play a key role in the proofs of Theorems 4 and 6 below.

Finally, let us mention that by using the approach of this section, one can, also, give an alternative proof of the formula \eqref{3}
of Erd\H{o}s and R\'{e}nyi \cite{E-R}.

\subsection{The limiting distribution of $D_{m,n}$ in the supercritical case}
\textbf{Theorem 3.} Let $\Delta_{m,n}$ be a random variable whose distribution function is given by \eqref{P3}, where $m = m(n)$, $n = 1, 2, \ldots$, is a sequence of positive integers such that
\begin{equation}
\frac{(\ln n)^3}{m(n)} \, \to 0
\label{AT1}
\end{equation}
(in other words, $m \gg \ln^3 n$).
Then, for any fixed $y \in \mathbb{R}$ we have
\begin{equation}
\lim_{n \rightarrow \infty} \mathbb{P}\left\{\frac{\Delta_{m,n} - n m - n \sqrt{m} \sqrt{2\ln n - \ln \ln n}}{n \sqrt{ m/2\ln n}} \leq y\right\}
= \exp \left(-\frac{e^{-y}}{2\sqrt{\pi}}\right)
\label{AT2}
\end{equation}
or, equivalently,
\begin{equation}
\sqrt{2 m \ln n} \left( \frac{\Delta_{m,n}}{n m} -1 \right)
- 2\ln n + \frac{\ln \ln n}{2} + \ln\left(2\sqrt{\pi}\right)
 \, \overset{d}{\longrightarrow} \, G,
\label{AT2a}
\end{equation}
where $G$ is the standard Gumbel random variable and the symbol $\overset{d}{\longrightarrow}$ denotes convergence in distribution.

\smallskip

\textit{Proof}. For typographical convenience we set
\begin{equation}
\lambda = \lambda(n) := m(n) + \sqrt{m(n)} \sqrt{2\ln n - \ln \ln n} + \frac{\sqrt{m(n)}}{\sqrt{2\ln n}} \, y,
\label{AT3}
\end{equation}
where $y$ is a fixed real number. Then, $\lambda > 0$ for all sufficiently large $n$ and \eqref{P3} implies
\begin{align}
\mathbb{P}\left\{\frac{\Delta_{m,n} - n m - n \sqrt{m} \sqrt{2\ln n - \ln \ln n}}{n \sqrt{ m/2\ln n}} \leq y\right\}
&= \mathbb{P}\left\{\Delta_{m,n} \leq n \lambda \right\}
\nonumber
\\
&= \left[ 1 - S_{m}(\lambda) e^{-\lambda} \right]^n,
\label{AT4}
\end{align}
where we have suppressed the dependence of $m$ and $\lambda$ on $n$ for typographical convenience.

We can, now, use \eqref{D6a} in order to compute the asymptotics of
$S_m(\lambda) e^{-\lambda}$. First we notice that from \eqref{AT3} we have $(\lambda - m) \to \infty$ and, furthermore,
\begin{equation}
\frac{\sqrt{m}}{\lambda - m} = \frac{1}{\sqrt{2\ln n - \ln \ln n} + \frac{y}{\sqrt{2\ln n}}}
\sim \frac{1}{\sqrt{2\ln n}},
\label{AT5}
\end{equation}
hence \eqref{D5a} is satisfied (for $x = \lambda$) and, consequently, formula \eqref{D6a} yields
\begin{equation}
S_m(\lambda) e^{-\lambda}
\sim \frac{1}{\sqrt{2\pi}} \left(\frac{\lambda}{m}\right)^m \frac{\sqrt{m}}{\lambda - m} \,
e^{-(\lambda-m)}.
\label{AT7}
\end{equation}
Next, by invoking \eqref{AT3}, formula \eqref{AT7} becomes
\begin{equation}
S_m(\lambda) e^{-\lambda}
\sim \left(1+ \frac{\sqrt{2\ln n - \ln \ln n} + \frac{y}{\sqrt{2\ln n}}}{\sqrt{m}}\right)^m
\frac{e^{-\sqrt{m} \left(\sqrt{2\ln n - \ln \ln n} + \frac{y}{\sqrt{2\ln n}}\right)}}{2\sqrt{\pi}\sqrt{\ln n}}.
\label{AT8}
\end{equation}
Now, in view of \eqref{AT1},
\begin{align}
&\left(1+ \frac{\sqrt{2\ln n - \ln \ln n} + \frac{y}{\sqrt{2\ln n}}}{\sqrt{m}}\right)^m
= e^{m \ln\left(1+ \frac{\sqrt{2\ln n - \ln \ln n} + \frac{y}{\sqrt{2\ln n}}}{\sqrt{m}}\right)}
 \nonumber
\\
&\qquad\qquad\qquad \sim e^{\sqrt{m} \left(\sqrt{2\ln n - \ln \ln n} + \frac{y}{\sqrt{2\ln n}}\right)
-\frac{1}{2} \left(\sqrt{2\ln n - \ln \ln n} + \frac{y}{\sqrt{2\ln n}}\right)^2}
 \nonumber
\\
&\qquad\qquad\qquad \sim e^{\sqrt{m} \left(\sqrt{2\ln n - \ln \ln n} + \frac{y}{\sqrt{2\ln n}}\right)
-\ln n + \frac{\ln \ln n}{2} - y}.
\label{AT9}
\end{align}
Substituting \eqref{AT9} in \eqref{AT8} we obtain
\begin{equation}
S_m(\lambda) e^{-\lambda}
\sim \frac{e^{-\ln n + \frac{\ln \ln n}{2} - y}}{2\sqrt{\pi}\sqrt{\ln n}}
= \frac{1}{n} \frac{e^{-y}}{2\sqrt{\pi}}
\label{AT10}
\end{equation}
and, therefore, formula \eqref{AT2} follows by using \eqref{AT10} in \eqref{AT4} and letting $n \to \infty$.
\hfill $\blacksquare$

\medskip

By a straightforward adaptation of the above proof we can cover the case where, instead of the assumption \eqref{AT1}, we allow the slightly more
general condition $m \gg (\ln n)^p$, with $1 < p < 3$. However, the formulas get considerably messier.

Theorem 3 is by itself interesting. However, our ultimate goal is to prove a similar statement for the variable $D_{m,n}$. In order to relate
$D_{m,n}$ to $\Delta_{m,n}$, we will use the following well-known lemma (see, e.g., \cite{C}, \cite{D}).

\medskip

\textbf{Converging Together Lemma.} Let $X_n$ and $Y_n$, $n = 1, 2, \ldots$, be two sequences of random variables such that
\begin{equation}
X_n \, \overset{d}{\longrightarrow} \, X
\qquad\quad \text{and} \qquad\quad
(Y_n - X_n)  \, \overset{d}{\longrightarrow} \, 0,
\label{DL1}
\end{equation}
where $X$ is some random variable. Then
\begin{equation}
Y_n \, \overset{d}{\longrightarrow} \, X.
\label{DL2}
\end{equation}
%
%

\medskip

We are now ready to give the limiting distribution of $D_{m,n}$.

\medskip

\textbf{Theorem 4.} Let $m = m(n)$, $n = 1, 2, \ldots$, be a sequence of positive integers such that $m \gg \ln^3 n$.
Then, for any fixed $y \in \mathbb{R}$ we have
\begin{equation}
\lim_{n \rightarrow \infty} \mathbb{P}\left\{\frac{D_{m,n} - n m - n \sqrt{m} \sqrt{2\ln n - \ln \ln n}}{n \sqrt{ m/2\ln n}} \leq y\right\}
= \exp \left(-\frac{e^{-y}}{2\sqrt{\pi}}\right)
\label{ATT2}
\end{equation}
or, equivalently,
\begin{equation}
\sqrt{2 m \ln n} \left( \frac{D_{m,n}}{n m} -1 \right)
- 2\ln n + \frac{\ln \ln n}{2} + \ln\left(2\sqrt{\pi}\right)
 \, \overset{d}{\longrightarrow} \, G,
\label{ATT2a}
\end{equation}
where $G$ is the standard Gumbel random variable.

\smallskip

\textit{Proof}. Let us set
\begin{equation}
Z_{m,n} := \frac{\Delta_{m,n} - D_{m,n}}{\sqrt{D_{m,n}}}.
\label{ATT3}
\end{equation}
Then
\begin{equation}
\mathbb{E}\left[Z_{m,n}^2\right]
= \mathbb{E}\left[\frac{\Delta_{m,n}^2}{D_{m,n}}\right] - 2\,\mathbb{E}\left[\Delta_{m,n}\right] + \mathbb{E}\left[D_{m,n}\right]
= \mathbb{E}\left[\frac{\Delta_{m,n}^2}{D_{m,n}}\right] - \mathbb{E}\left[D_{m,n}\right],
\label{ATT4}
\end{equation}
where the second equality follows from \eqref{5}. Now, in view of \eqref{A9} and \eqref{P5},
\begin{equation}
\mathbb{E}\left[\left.\frac{\Delta_{m,n}^2}{D_{m,n}} \, \right|  D_{m,n}\right]
= \frac{1}{D_{m,n}} \, \mathbb{E}\left[\Delta_{m,n}^2 \, | \,  D_{m,n}\right]
= \frac{D_{m,n}^{(2)}}{D_{m,n}} = D_{m,n} + 1
\label{ATT5}
\end{equation}
and hence
\begin{equation}
\mathbb{E}\left[\frac{\Delta_{m,n}^2}{D_{m,n}}\right] = \mathbb{E}\left[D_{m,n}\right] + 1.
\label{ATT6}
\end{equation}
Therefore, by substituting \eqref{ATT6} in \eqref{ATT4} we obtain that
\begin{equation}
\mathbb{E}\left[Z_{m,n}^2\right] = 1
\label{ATT7}
\end{equation}
(actually, in view of \eqref{P4} and the fact that $D_{m,n} \geq mn$, one can show that $Z_{m,n}$ converges in distribution to a standard normal
random variable, as $m \to \infty$ or $n \to \infty$).

Now, from \eqref{ATT3} we get
\begin{equation*}
\left|\Delta_{m,n} - D_{m,n}\right| = \left|Z_{m,n}\right| \sqrt{D_{m,n}}
\end{equation*}
and, hence, in view of \eqref{ATT7},
\begin{equation}
\mathbb{E}\big[\left|\Delta_{m,n} - D_{m,n}\right|\big] = \mathbb{E}\left[\left|Z_{m,n}\right| \sqrt{D_{m,n}}\right]
\leq \sqrt{\mathbb{E}\left[Z_{m,n}^2\right] \mathbb{E}\left[D_{m,n}\right]} = \sqrt{\mathbb{E}\left[D_{m,n}\right]}.
\label{ATT8}
\end{equation}
Thus, in view of \eqref{M25},
\begin{equation}
\mathbb{E}\left[\frac{\left|\Delta_{m,n} - D_{m,n}\right|}{n \sqrt{ m/2\ln n}}\right]
\leq \frac{\sqrt{\mathbb{E}\left[D_{m,n}\right]}}{n \sqrt{ m/2\ln n}} \to 0.
\label{ATT9}
\end{equation}
It follows that we can apply the Converging Together Lemma to the sequences
\begin{equation*}
X_n := \frac{\Delta_{m,n} - n m - n \sqrt{m} \sqrt{2\ln n - \ln \ln n}}{n \sqrt{m/2\ln n}}
\end{equation*}
and
\begin{equation*}
Y_n := \frac{D_{m,n} - n m - n \sqrt{m} \sqrt{2\ln n - \ln \ln n}}{n \sqrt{ m/2\ln n}}
\end{equation*}
and conclude that their limiting distributions coincide.
\hfill $\blacksquare$

\medskip

\textbf{Remark 4.} In the somehow extreme case where $m \to \infty$, while $n$ stays fixed, the limiting distribution of $\Delta_{m,n}$
follows directly from \eqref{P1}: Since each $T_j$, $j = 1, 2, \ldots, n$ is Erlang with parameters $m$ and $1/n$, the Central Limit Theorem
yields
\begin{equation}
\frac{T_j - mn}{n\sqrt{m}}  \, \overset{d}{\longrightarrow} \, Z,
\qquad
m \to \infty,
\label{ATT11}
\end{equation}
where $Z$ is the standard normal random variable. Then, formula \eqref{P1} and the independence of the $T_j$'s imply immediately that
\begin{equation}
\frac{\Delta_{m,n} - n m}{n \sqrt{m}}
 \, \overset{d}{\longrightarrow} \, \max\{Z_1, Z_2, \ldots, Z_n\},
 \qquad
m \to \infty,
\label{ATT12}
\end{equation}
where $Z_1, Z_2, \ldots, Z_n$ are independent standard normal variables (also, in view of \eqref{5}, we can obtain that
$\mathbb{E}[D_{m,n}] = \mathbb{E}\left[ \Delta_{m,n} \right] \sim nm$ as $m \to \infty$). However, in the case of a fixed $n$, formula \eqref{ATT9}
in the proof of Theorem 4 fails and, hence, we cannot conclude that $D_{m,n}$ and $\Delta_{m,n}$ have the same limiting distributions. Actually,
in the trivial case $n=1$ we obviously have $D_{m,1} = m$, while from \eqref{ATT12} we see that
$(\Delta_{m,1} - m)/\sqrt{m} \, \overset{d}{\longrightarrow} \, Z$.

\subsection{The limiting distribution of $D_{m,n}$ in the critical case}
We will now consider the case
\begin{equation}
m = m(n) = \beta \ln n + b(n) \ln n,
\qquad \text{where }\;
b(n) = o\left(\frac{1}{\sqrt{\ln n}}\right),
\label{CC1}
\end{equation}
i.e. $m = \beta \ln n + o\big(\sqrt{\ln n}\,\big)$ (as usual, $\beta > 0$ is a fixed constant). The restriction on the order of $b(n)$ is imposed in order to
keep the formulas relatively simple (we believe that by a straightforward adaptation of the proof of Theorem 5 below we can cover the more general
case $b(n) \ll (\ln n)^{-p}$ for a fixed $p \in (0, 1/2)$; however, the formulas will get considerably messier).

\medskip

\textbf{Theorem 5.} Let $\Delta_{m,n}$ be a random variable whose distribution function is given by \eqref{P3}, where $m = m(n)$, $n = 1, 2, \ldots$, is a sequence of positive integers satisfying \eqref{CC1}. Then, for any fixed $y \in \mathbb{R}$ we have
\begin{align}
\lim_{n \rightarrow \infty}
&\mathbb{P}\left\{\frac{\Delta_{m,n} - \alpha n \ln n - \frac{\alpha(\alpha-\beta-1)}{\beta(\alpha-\beta)} \, b(n) n \ln n + \frac{\alpha}{2(\alpha-\beta)} \, n \ln \ln n}{\frac{\alpha}{\alpha-\beta} \, n}
\leq y\right\}
\nonumber
\\
&= \exp \left(-\frac{1}{\sqrt{2\pi}} \, \frac{\sqrt{\beta}}{\alpha - \beta} \, e^{-y}\right),
\label{T5AT2}
\end{align}
where $\alpha$ is given by \eqref{CL2} (recall that $\alpha > \beta$).

\smallskip

\textit{Proof}. We will follow the approach of the proof of Theorem 3.

For typographical convenience we set
\begin{equation}
\chi = \chi(n)
:= \alpha \ln n + \frac{\alpha(\alpha-\beta-1)}{\beta(\alpha-\beta)} \, b(n) \ln n
- \frac{\alpha}{2(\alpha-\beta)} \ln \ln n + \frac{\alpha}{\alpha-\beta} \, y,
\label{T5AT3}
\end{equation}
where $y$ is a fixed real number. Since $b(n) = o(1)$, the quantity $\chi$ is positive for all sufficiently large $n$ and, hence, \eqref{P3} implies
\begin{align}
&\mathbb{P}\left\{\frac{\Delta_{m,n} - \alpha n \ln n - \frac{\alpha(\alpha-\beta-1)}{\beta(\alpha-\beta)} \, b(n) n \ln n + \frac{\alpha}{2(\alpha-\beta)} \, n \ln \ln n}{\frac{\alpha}{\alpha-\beta} \, n}
\leq y\right\}
\nonumber
\\
&= \mathbb{P}\left\{\Delta_{m,n} \leq n \chi \right\} = \left[ 1 - S_{m}(\chi) e^{-\chi} \right]^n.
\label{T5AT4}
\end{align}
Next, we write \eqref{T5AT3} as
\begin{equation}
\chi(n) = \alpha \ln n + a(n) \ln n,
\label{T5AT3a}
\end{equation}
where
\begin{equation}
a = a(n) := \frac{\alpha(\alpha-\beta-1)}{\beta(\alpha-\beta)} \, b(n)
- \frac{\alpha}{2(\alpha-\beta)} \frac{\ln \ln n}{\ln n} + \frac{\alpha}{\alpha-\beta} \frac{y}{\ln n},
\label{T5AT3b}
\end{equation}
so that, our assumption \eqref{CC1} implies
\begin{equation}
a(n) = o\left(\frac{1}{\sqrt{\ln n}}\right)
\label{T5AT3c}
\end{equation}
(in particular, $a(n) \to 0$ as $n \to \infty$). We can, therefore, invoke formulas \eqref{CL8} and \eqref{CL8a} appearing in the proof of Lemma 1
and obtain
\begin{equation}
S_m(\chi) e^{-\chi} \sim \frac{1}{n} \, \frac{1}{\sqrt{2\pi}} \, \frac{\sqrt{\beta}}{\alpha - \beta} \, e^{A(n) \ln n},
\label{T5CL8}
\end{equation}
where
\begin{align}
A(n) &:=(\ln\alpha - \ln\beta + 1) b - (\beta + b) \ln\left(1 + \frac{b}{\beta}\right)
\nonumber
\\
&- a + (\beta + b) \ln\left(1 + \frac{a}{\alpha}\right) - \frac{1}{2} \frac{\ln\ln n}{\ln n}
\label{T5CL8a}
\end{align}
(here, $a = a(n)$ is, of course, given by \eqref{T5AT3b}).

Now, in view of \eqref{CC1} and \eqref{T5AT3c}, formula \eqref{T5CL8a} implies
\begin{equation}
A(n) = (\ln\alpha - \ln\beta)b - \frac{1}{2} \frac{\ln\ln n}{\ln n} - \frac{\alpha-\beta}{\alpha} \, a + o\left(\frac{1}{\ln n}\right).
\label{T5CL9}
\end{equation}
We can, then, substitute \eqref{T5AT3b} in \eqref{T5CL9} and get
\begin{equation}
A(n) = (\ln\alpha - \ln\beta)b - \frac{\alpha-\beta-1}{\beta} \, b - \frac{y}{\ln n} + o\left(\frac{1}{\ln n}\right),
\label{T5CL10}
\end{equation}
or, in view of \eqref{CL2},
\begin{equation}
A(n) = - \frac{y}{\ln n} + o\left(\frac{1}{\ln n}\right),
\label{T5CL11}
\end{equation}
Thus, by substituting \eqref{T5CL11} in \eqref{T5CL8} we obtain
\begin{equation}
S_m(\chi) e^{-\chi} \sim \frac{1}{n} \, \frac{1}{\sqrt{2\pi}} \, \frac{\sqrt{\beta}}{\alpha - \beta} \, e^{-y}.
\label{T5AT12}
\end{equation}
Therefore, formula \eqref{T5AT2} follows by using \eqref{T5AT12} in \eqref{T5AT4} and then letting $n \to \infty$.
\hfill $\blacksquare$

\medskip

Let us notice that formula \eqref{T5AT2} can be, also, expressed equivalently as
\begin{equation}
\frac{\alpha-\beta}{\alpha} \, \frac{\Delta_{m,n}}{n} - (\alpha-\beta) \ln n - \frac{\alpha-\beta-1}{\beta} \, b(n) \ln n
+ \frac{\ln \ln n}{2} + C
 \, \overset{d}{\longrightarrow} \, G,
\label{T5AT2a}
\end{equation}
with
\begin{equation}
C := \frac{1}{2}\ln\left(\frac{2\pi(\alpha-\beta)^2}{\beta}\right)
\label{T5AT2b}
\end{equation}
(as usual, $G$ is the standard Gumbel random variable). Furthermore, in view of \eqref{CC1}, formula \eqref{T5AT2a} can be also written as
\begin{equation}
\frac{\alpha-\beta}{\alpha} \, m \left(\frac{\Delta_{m,n}}{n m} - 1\right) + \frac{b(n)}{\beta} \, \ln n
+ \frac{\ln \ln n}{2} + C
 \, \overset{d}{\longrightarrow} \, G.
\label{T5AT2c}
\end{equation}

As in the supercritical case, in the critical case too, under the above normalization $\Delta_{m,n}$ and $D_{m,n}$ have the same limiting
distributions:

\medskip

\textbf{Theorem 6.} Let $m = m(n)$, $n = 1, 2, \ldots$, be a sequence of positive integers satisfying \eqref{CC1}. Then, for any fixed $y \in \mathbb{R}$ we have
\begin{align}
\lim_{n \rightarrow \infty}
&\mathbb{P}\left\{\frac{D_{m,n} - \alpha n \ln n - \frac{\alpha(\alpha-\beta-1)}{\beta(\alpha-\beta)} \, b(n) n \ln n + \frac{\alpha}{2(\alpha-\beta)} \, n \ln \ln n}{\frac{\alpha}{\alpha-\beta} \, n}
\leq y\right\}
\nonumber
\\
&= \exp \left(-\frac{1}{\sqrt{2\pi}} \, \frac{\sqrt{\beta}}{\alpha - \beta} \, e^{-y}\right),
\label{T5ATc}
\end{align}
where $\alpha$ is given by \eqref{CL2}. Equivalently,
\begin{equation}
\frac{\alpha-\beta}{\alpha} \, m \left(\frac{D_{m,n}}{n m} - 1\right) + \frac{b(n)}{\beta} \, \ln n
+ \frac{\ln \ln n}{2} + C
 \, \overset{d}{\longrightarrow} \, G,
\label{T5ATd}
\end{equation}
where $C$ is given by \eqref{T5AT2b} and $G$ is the standard Gumbel random variable.
\medskip

The proof of Theorem 6 is exactly the same as the proof of Theorem 4. The only difference here is that, in order to justify the analog of estimate
\eqref{ATT9}, instead of \eqref{M25}, we now use \eqref{MM25}.

\subsection{Comparison of the formulas \eqref{ATT2a} and \eqref{T5ATd}}
As we saw in the previous subsections, in the critical case $m = \beta \ln n + o\big(\sqrt{\ln n}\big)$ the limiting behavior of $D_{m,n}$ is given
by \eqref{T5ATd}, while in the supercritical case $m \gg \ln^3 n $ the limiting behavior of $D_{m,n}$ is given by \eqref{ATT2a}. In this short
subsection we will enquire whether these two behaviors can, in some sense, be ``bridged."

Heuristically, one expects that the transition from the critical to the supercritical case can be observed by letting $\beta \to \infty$.

Let us first notice that the left-hand side of \eqref{ATT2a} as well as the left-hand side of \eqref{T5ATd} consist of four terms. Of these four
terms, the third, namely $\ln \ln n/2$, is common in both formulas, while the second terms do not agree, although they both contain the
factor $\ln n$.

From the fact that $\alpha > \beta$ it follows that if $\beta \to \infty$, then $\alpha \to \infty$. Now formula \eqref{CL2} implies
\begin{equation}
\frac{\alpha}{\beta} - 1 = \ln\left(\frac{\alpha}{\beta}\right) + \frac{1}{\beta},
\label{E1}
\end{equation}
hence
\begin{equation}
\frac{\alpha}{\beta} - 1 = \ln\left(\frac{\alpha}{\beta}\right) + o(1)
\qquad\text{as }\;
\beta \to \infty
\label{E2}
\end{equation}
from which we get that
\begin{equation}
\frac{\alpha}{\beta} \to 1
\qquad\text{as }\;
\beta \to \infty.
\label{E3}
\end{equation}
Let us set
\begin{equation}
u := \frac{\alpha}{\beta} - 1.
\label{E4}
\end{equation}
Then, \eqref{E1} can be written as
\begin{equation}
u = \ln\left(1+u\right) + \frac{1}{\beta}.
\label{E5}
\end{equation}
Since (in view of \eqref{E3}) $u = o(1)$ as $\beta \to \infty$, formula \eqref{E5} implies
\begin{equation*}
u = u - \frac{u^2}{2} + \frac{1}{\beta} + O\left(u^3\right),
\qquad
\beta \to \infty
\end{equation*}
i.e.
\begin{equation}
u = \sqrt{\frac{2}{\beta}} + O\left(u^{3/2}\right),
\qquad
\beta \to \infty.
\label{E6}
\end{equation}
Since $u \gg u^{3/2}$ as $\beta \to \infty$, formula \eqref{E6} can be written as
\begin{equation}
u = \sqrt{\frac{2}{\beta}} + O\left(\frac{1}{\beta^{3/4}}\right),
\qquad
\beta \to \infty,
\label{E7}
\end{equation}
or, in view of \eqref{E5},
\begin{equation}
\frac{\alpha}{\beta} - 1 = \sqrt{\frac{2}{\beta}} + O\left(\frac{1}{\beta^{3/4}}\right),
\qquad
\beta \to \infty,
\label{E8}
\end{equation}
which implies immediately that
\begin{equation}
\frac{\alpha - \beta}{\sqrt{\beta}} = \sqrt{2} + O\left(\frac{1}{\beta^{1/4}}\right),
\qquad
\beta \to \infty.
\label{E9}
\end{equation}
Using \eqref{E9} in \eqref{T5AT2b} we obtain that
\begin{equation}
C = \frac{1}{2}\ln\left(\frac{2\pi(\alpha-\beta)^2}{\beta}\right) \to \frac{1}{2}\ln\left(4\pi\right) = \ln\left(2\sqrt{\pi}\right)
\qquad\text{as }\;
\beta \to \infty.
\label{E10}
\end{equation}
Therefore, the fourth term of the left-hand side of \eqref{T5ATd} (namely $C$) approaches the corresponding (fourth) term of the left-hand side of
\eqref{ATT2a}, namely $\ln\left(2\sqrt{\pi}\right)$, as $\beta \to \infty$.

As for the first term of the left-hand side of \eqref{T5ATd} we have (in view of \eqref{E3} and \eqref{E9})
\begin{equation}
\frac{\alpha-\beta}{\alpha} \, m
\sim \frac{\alpha-\beta}{\alpha} \, \beta \ln n
\sim (\alpha-\beta) \ln n
= \frac{\alpha-\beta}{\sqrt{\beta}} \, \sqrt{\beta} \ln n
\sim \sqrt{2\beta} \ln n
\label{E11}
\end{equation}
and this is in ``asymptotical agreement" with the corresponding factor of the first term of the left-hand side of \eqref{ATT2a} (in the case
where $m \sim \beta \ln n$) since
\begin{equation}
\sqrt{2m \ln n}
\sim \sqrt{2\beta \ln^2 n}
\sim \sqrt{2\beta} \ln n.
\label{E12}
\end{equation}

\textbf{Acknowledgments.} The authors wish to thank an anonymous referee for informing them about G.I. Ivchenko's papers \cite{I1} and \cite{I2}.

\end{document}